\documentclass[12pt,a4paper]{article}

\usepackage{amsmath}
\usepackage{amsthm}
\usepackage[ruled,vlined,linesnumbered,inoutnumbered]{algorithm2e}
\usepackage{bm}
\usepackage[caption=false,font=footnotesize]{subfig}
\usepackage{floatflt}
\usepackage{epstopdf}
\usepackage{graphicx}

\title{On enhancing efficiency and accuracy of particle swarm optimization algorithms}
\date{}
\begin{document}

\maketitle

\centerline{\scshape  Silvano Chiaradonna, Felicita Di Giandomenico and Nadir Murru}
 \medskip
{\footnotesize
\centerline{ISTI -- CNR} 
\centerline{Via G. Moruzzi 1, 56124, Pisa, ITALY}
\centerline{silvano.chiaradonna@isti.cnr.it; f.digiandomenico@isti.cnr.it; nadir.murru@isti.cnr.it} }
\medskip

%%\centerline{(Communicated by xxxx)}

\medskip

\begin{abstract} 
{\em The particle swarm optimization (PSO) algorithm has been recently introduced in the non--linear programming, becoming widely studied and used in a variety of applications. Starting from its original formulation, many variants for improvement and specialization of the PSO have been already proposed, but without any definitive result, thus research in this area is nowadays still rather active. This paper goes in this direction, by proposing some modifications to the basic PSO algorithm, aiming at enhancements in aspects that impact on the efficiency and accuracy of the optimization algorithm. In particular, variants of PSO based on fuzzy logics and Bayesian theory have been developed, which show better, or competitive, performances when compared to both the basic PSO formulation and a few other optimization algorithms taken from the literature.}
\end{abstract}

\section{Introduction}

Many techniques, mainly based on metaheuristic algorithms, have been developed for solving non--linear programming. In particular, chaotic algorithms \cite{Sh}, evolutionary programming \cite{Wu, Yu}, genetic algorithms \cite{Zhang, Bak}, tabu search algorithm \cite{Ku, Lan}, have been widely used in these studies. Recently, the particle swarm optimization (PSO) algorithm \cite{Kennedy} has been introduced in the non--linear programming becoming widely studied and used in a variety of applications \cite{Zhao, Esmin, KZAW13}.

The PSO algorithm performs a metaheuristic search based on competition and cooperation among particles (which represent the search variables) belonging to an initial swarm or population that covers the search space. After the basic technique was defined, numerous proposals for improvement and specialization of the PSO have appeared (refer \cite{Poli2008} for a partial survey). Still, further investigations are nowadays of interest and solicited in PSO, to better respond to application needs in a growing variety of sectors, with related emphasis on different aspects, such as efficiency and accuracy.

In this paper we propose some modifications to the basic PSO algorithm, aiming at enhancements in the following three aspects that determine the efficiency and accuracy of the optimization algorithm: 
\begin{itemize}
\item  the setting of the initial population, based on the partitioning of the search space in disjoint intervals, one for each particle, such that the initial position for each different particle is selected from a different interval,
\item the setting of the acceleration constants and of the inertia weight parameter, defined by using a fuzzy logics based strategy,
\item the calculation of the current position of each particle, based on a Bayesian approach.
\end{itemize}

By adopting approaches very suited to support the dynamics of the PSO algorithm, such as the fuzzy logics and the Bayesian theory, this paper offers novel solutions which improve on the basic PSO formulation, thus potentially resulting more appealing to be employed, as suggested by the evaluation studies already performed. Although the literature on this family of optimazation algorithms is very abundant, as already recalled, our goal is to enrich the set of alternatives to enhance the choice of the best solution that fits specific applications needs and characteristics.

The rest of the paper is structured as follows. A brief recall of the basic PSO algorithm is presented in Section~\ref{sec:psobasic}. In Section~\ref{sec:innovativeaspects}, the motivation and innovative aspects of the three proposed modifications to the basic PSO algorithm are discussed. In the next three sections, we develop our enhancements to the basic PSO algorithm. In particular, the setting proposed for the initial population is described in Section~\ref{sec:sample}; the approach proposed to set the acceleration constants and the inertia weight is described in Section~\ref{sec:fuzzy}, while Section~\ref{sec:bayes} contains the new method to move the control variables in the search space. By composing the proposed improvements, two variants of the PSO algorithm are assembled in Section~\ref{sec:simulations}, which are then evaluated and compared both with the basic PSO algorithm in terms of indicators representative of their accuracy and efficiency, and with other solutions taken from the literature when employed to support power flow optimization in a simple but realistic case study in the electrical power sector. Finally, Section ~\ref{sec:conclusions} is devoted to the conclusions.

\section{Overview of the basic PSO algorithm} \label{sec:psobasic}

In this section, the basic PSO algorithm is briefly described, to provide the context for the variants that will be presented in the next sections.
In general terms, an optimization problem is the problem of minimizing (or maximizing) an objective function $f(x_1 \dots, x_d, y_1, \dots, y_s)$, where $x_i$ are the search variables, representing the parameters of the problem that are to be optimized, and $y_i$ are the state variables, representing the parameters that can be derived as  function of the search variables. The variables $x_i$ and $y_i$ can be bounded and related together by means of constraint equations or disequations. 
For the sake of simplicity, in the following, we omit variables $y_i$ from the notation of the objective function and we only consider search variables $x_i$.

The PSO algorithm is a metaheuristic global optimization method which tries to find those values of the search variables $x_1,...x_d$ for which the value of the cost function $f(x_1,...x_d)$ is minimum (or maximum). The search is based on a discrete-time stochastic process describing the competition and cooperation between a population or swarm of particles that move in a search-variables space. 
The PSO algorithm considers a fixed size swarm of $D$ particles that randomly change position with velocity $\mathbf{v}_p(k)$ at each discrete instant of time (step) $k$, according to their own experience and that of their neighbors. 
The position of each particle $p$ at step $k$ is a vector $\mathbf{x}_p(k)=(x_{p1}(k),..., x_{pd}(k))$.
The solution of the optimization problem found by PSO is the position $\mathbf{b}$, reached at a certain step by of a particle for the swarm, for which the value of the cost function is minimum, i.e., $f(\mathbf{b})=\min\{ f(\mathbf{x}_p(k)) | p=1, \dots, D, k=1, 2, \dots\}$.
Algorithm~\ref{alg:BasicPSOalgorithm} describes formally in pseudo--code the basic PSO procedure of finding 
the optimal vector $\mathbf{b}$, as extracted from \cite{Kennedy}.
\begin{algorithm}[!htb]
\DontPrintSemicolon
\KwData{$d,D,w,c_1,c_2,k_{max}$ and bounds $\mathbf{x}^{min},\mathbf{x}^{max}$\;
Search variables $\mathbf{x}=(x_1,...x_d)$ and cost function $f(\mathbf{x})$\;
}
\KwResult{
  $\mathbf{b}$, $f(\mathbf{b})$
}
\For(\tcp*[f]{for each particle}){$p=1,\dots,D$}{$\mathbf{v}_p(0) \gets \mathbf{0}$ \nllabel{algl:startingv}  \tcp*{initial velocity} $\mathbf{x}_p(0) \gets Rnd(\mathbf{x}^{min},\mathbf{x}^{max})$ \nllabel{algl:startingx} \tcp*{initial random position} $\mathbf{l}_p,\mathbf{b} \gets \mathbf{x}_p(0)$ } 
\Begin{
\For(\tcp*[f]{for each step}){$k=1, \ldots, k_{max}$}{ 
  \For(\tcp*[f]{for each particle}){$p=1,\dots,D$}{
    \If{$f(\mathbf{x}_p(k-1)) < f(\mathbf{l}_p)$}{
      $\mathbf{l}_p \gets \mathbf{x}_p(k-1)$ \tcp*{new personal best position} 
      \nllabel{algl:l}
    }
    \If{$f(\mathbf{l}_p)<f(\mathbf{b})$}{
      $\mathbf{b} \gets \mathbf{l}_p$\tcp*{new global best position} }  \nllabel{algl:b}    
  }
  \For(\tcp*[f]{for each particle}){$p=1,\dots,D$}{
    $r_1 \gets \mathcal U(0,1)$ \tcp*{new uniform random value}\nllabel{algl:r1}
    $r_2 \gets \mathcal U(0,1)$ \tcp*{new uniform random value}
    $\mathbf{v}_p(k) = w \mathbf{v}_p(k-1) + c_1 r_1 (\mathbf{l}_p - \mathbf{x}_p(k-1)) + c_2 r_2 (\mathbf{b} - \mathbf{x}_p(k-1) )$\; \nllabel{algl:vk}
    $\label{part} \mathbf{x}_p(k) = \mathbf{x}_p(k-1) + \mathbf{v}_p(k)$\;   \nllabel{algl:xk}
    \For(\tcp*[f]{for each search variable}){$i=1, \ldots, d$}{
      \If{$x_{pi}(k) < x_i^{min}$}{$x_{pi}(k)\gets x_i^{min} $}
      \ElseIf{$x_{pi}(k) > x_i^{max}$}{$x_{pi}(k)\gets x_i^{max} $}
    }
  } 
}
}
\caption{Basic PSO algorithm}
\label{alg:BasicPSOalgorithm}
\end{algorithm}

The starting velocity and position of the particles at step $0$ are initialized at lines~\ref{algl:startingv} and~\ref{algl:startingx}, respectively.
The personal best position $\mathbf{l}_p$ for the particle $p$ at step $k$ is the best position the particle $p$ has visited since the first step, i.e., 
the position that minimizes the objective function among all the positions of $p$ for different values of $k$. It is updated at line~\ref{algl:l}.
The global best position $\mathbf{b}$ at each step $k$ is the best position discovered by any of the particles so far. It is usually calculated as the personal best position and it is obtained as the personal best position at step $k$ that minimizes the objective function, as shown at line~\ref{algl:b}, where $\mathbf{b}$ is updated.

The position $\mathbf{x}_p(k)$ of the particle $p$ at step $k$ is updated at line~\ref{algl:xk}, based on the random value assigned to the velocity $\textbf v_p(k)$, as shown from line~\ref{algl:r1} to line~\ref{algl:vk}. 
The velocity vector drives the optimization process and results from the sum of three different components, as shown at line~\ref{algl:vk}: the momentum component, which is the previous velocity, 
the cognitive component, which is proportional to the distance of the particle from the best position it has ever visited, and, finally, the social component which is proportional to the particle's distance from the best position where any of the swarm's particles has ever been.

The PSO algorithm has some critical points that heavily influence its performances, like the initial position of the particles and the choice of some parameters, such as inertia weight $w$ and the acceleration constants $c_1$ and $c_2$, used to scale the contribution of the moment to the velocity, and of the cognitive and social components, respectively~\cite{E07}. Indeed, 
if the initial position of the particles is a well-distributed cover of the search space, then there is a higher probability to avoid local best position and to approach the global optimum. Moreover, parameters $w, c_1, c_2$ heavily influence the velocity of the particles, i.e., their movements in the search space. 

In the next sections we present and evaluate three variants of the PSO algorithm, aiming at enhancements in the aspects just recalled that determine the efficiency and accuracy of the optimization algorithm.

\section{Three new variants of the PSO algorithm: motivation and innovative aspects}\label{sec:innovativeaspects}

In this section, we discuss the motivation and the innovative aspects of the three proposed variants for the PSO algorithm.

The efficiency of PSO is influenced by the initial position of the particles, i.e. by how well particles are distributed over the search space, since if regions of the search space are not covered by the initial positions, the PSO will have difficulty in finding the optimum if it is located within an uncovered region~\cite{E07}. Usually, the initial position of each particle is randomly generated, with uniform distribution, from the set of all possible positions of the particles. In this case, it is important to use a good pseudo--random generator, otherwise the performances heavily degrade. In \cite{Mu}, the authors propose to reinitialize the particle that has the worst performance (i.e., corresponding to the biggest objective function value). In \cite{Pani}, re--initialization of a certain percentage of the total population is considered if after some steps there is no significative improvement. In \cite{Zhao}, an orthogonal design is proposed in order to sample the initial positions. In \cite{WuZhang}, the authors determine the initial position of the particles by means of the Tent chaotic map. In Section~\ref{sec:sample}, we propose an original algorithm that forces a uniform distribution for the position of the initial particles, dividing the search space into small intervals where we randomly generate the initial positions.

As discussed in the previous section, the values of parameters $w, c_1, c_2$ affect the performance of the PSO algorithm. Standard values for the parameters of the PSO algorithm are $c_1=c_2=1.496172$ and $w=0.72984$ \cite{Bratt}. However, in many works (like \cite{Bansal}, \cite{Jiang}, \cite{Nick} and  \cite{Pani}) values of $c_1, c_2, w$ are not constant, but they change during the execution of the algorithm. This is mainly due to the fact that at the initial steps of the algorithm large values for $w$ are preferable so that particles have a greater possibility of movement \cite{Ojha}. On the contrary, small values for $w$ are preferable towards the final steps of the algorithm. 
Moreover, other different components could be chosen for modifying the values of $w, c_1, c_2$, like, e.g., the distance between the current particles and the best particles. Since, there are no defined rules to determine the relationship between parameters $w, c_1, c_2$ and other components (such as the current step of the algorithm), in  Section~\ref{sec:fuzzy}, we use fuzzy logic in order to manage these relationships. Indeed, fuzzy logics is congenial to capture and to code expert--based knowledge in view of performing targeted simulations. 
Usually the fuzzy logic--based systems are tuned using heuristic criteria (see, e.g., \cite{Der, Silva}). Fuzzy logic has been already used in various ways for improving the PSO algorithm~\cite{Abdel, Niknam}. In~\cite{SE01}, a fuzzy system is defined to dynamically adapt the inertia weight of the PSO. It is based on two input variables,  which measure 
the current best performance and the current inertia weight, and on an output variable, which is the change of the inertia weight. Our approach differs from previous usages of fuzzy logics (as in \cite{Niknam, SE01}) in the choice of the input fuzzy variables, in the tuning of the membership functions and in the fuzzy rules that combine input/output fuzzy variables together with their membership functions. 

Finally, in  Section~\ref{sec:bayes} we propose an innovative method to update the position of the particles at each step of the algorithm. The core idea of the basic PSO algorithm is to move the particles around the search space using information provided by personal and global best positions as shown in Algorithm~\ref{alg:BasicPSOalgorithm}. Here, we propose to move the particles in the search space by balancing the current position of each particle with the current personal and global best positions by means of equations based on the Kalman filter \cite{Kalman}. 
The Kalman filter uses a Bayesian approach providing a posterior estimation for a system's state from the current measurement and the prior estimation of that system's state by means of the Bayes theorem. 
In typical applications of the Kalman filter, it takes a sequence of measurements over time, containing  noise (random variations) and other inaccuracies, and produces statistically optimal estimates of unknown variables.
In our context, such a Bayesian approach has been reproduced considering the current position of the particles and the current global best position as
the current measurement and the prior estimation, respectively. The posterior estimation is the new global best position of the particles that should be an improved estimation for the position of the unknown global best.

\section{Variant 1: focus on the initial position of the particles} \label{sec:sample}

Here, we propose an easy sampling method that guarantees a good uniform distribution of the initial position of the particles in the search space, thus improving the use of pseudo--random generators.  

It is based on the partitioning of each dimension of the search space in disjoint intervals, one for each particle, such that the initial position of each different particle is randomly selected for each dimension from a different random interval.

Let $x^{min}_i$ and $x^{max}_i$ be the lower and upper bound of each search variable, i.e.,
\[x^{min}_i\leq x_{pi}(k) \leq x^{max}_i, \quad \forall p=1,\dots,D,\quad \forall i=1,\dots,d,\quad \forall k.\]
Each interval $[x^{min}_i,x^{max}_i]$ is partitioned into $D$ disjoint subintervals, one for each particle $p$
\[\left[ x^{min}_i+(p-1)\frac{x^{max}_i-x^{min}_i}{D}, x^{min}_i+p\frac{x^{max}_i-x^{min}_i}{D} \right], \quad \forall p=1,\dots,D.\]
\begin{algorithm}[!htb]
\DontPrintSemicolon
\KwData{$d, D, x^{min}_i, x^{max}_i$, for $i=1,\dots, d$}
\KwResult{
  $x_{pi}(0)$ for $p=1, \ldots ,D$ and $i=1,...,d$
}
\Begin{
\For(\tcp*[f]{for each search variable}){$i=1,\ldots,d$}{
  $\rho$ $\gets$ random permutation of $(1,\ldots,D)$\; \nllabel{algp:perm}
  \For(\tcp*[f]{for each particle}){$p=1,\ldots,D$}{
    $x_{pi}(0) \gets \mathcal U\left(x^{min}_i+(\rho(p)-1)\frac{x^{max}_i-x^{min}_i}{D}, x^{min}_i+\rho(p)\frac{x^{max}_i-x^{min}_i}{D}\right)$ \nllabel{algp:xpi}}}

}
\caption{Initial position of the particles for the PSO algorithm \label{alg:PSOinitial}}
\end{algorithm}
The initial position $x_{pi}(0)$ of each different particle $p$ for each dimension $i$ is uniformly randomly selected from a randomly selected different subinterval, as shown in Algorithm~\ref{alg:PSOinitial}.
At line~\ref{algp:xpi} it is shown that the position of the particle $p$ is selected by the $\rho(p)$--th subinterval, where $\rho()$ is the random permutation function of the indexes of the particles, as shown at line~\ref{algp:perm}.
Using this approach, the location of each different particle is selected for each dimension (search variable) from different partitions of the range of values for that dimension, improving the distribution of the particles in the search space.
This approach also protects against weaknesses of bad pseudo--random generators, which might result in a poor cover of the search space.

\section{Variant 2: fuzzy logic--based inertia weight and acceleration constants} \label{sec:fuzzy}

Usually, as shown in \cite{Bratt}, the following relation among $w, c_1, c_2$ holds:
\[w=\frac{2}{\lvert 2-\phi-\sqrt{\phi^2-4\phi} \rvert},\quad \phi=w(c_1+c_2),\]
from which, using equal values for $c_1, c_2$, we obtain
\begin{equation} \label{c1c2} c_1=c_2=\frac{(w+1)^2}{2}.  \end{equation}
Some works, like \cite{Nick}, propose a varying value for $w$:
\[w=\frac{k_{max}-k}{k_{max}}(w_{max}-w_{min})+w_{min},\]
where $k$ is the current step of the algorithm, $k_{max}$ is the maximum number of steps, and usually $w_{min}=0.8$ and $w_{max}=1.2$. In \cite{Pani} and \cite{Zheng} different formulas for $w$ have been proposed and in \cite{Jiang} similar formulas have been used for the coefficients $c_1$ and $c_2$. In \cite{Bansal}, a good survey on the problem of determining the parameters $w, c_1, c_2$ is provided. 

Here, we propose a novel fuzzy strategy to determine the parameters $w, c_1, c_2$ based on the steps number and the current particle position. 
Generally, a fuzzy strategy is composed by
\begin{itemize}
\item fuzzy variables (input and output), fuzzy sets and membership functions, 
\item fuzzy rules that relate input and output variables,
\item a fuzzy inference engine that combines the fuzzy rules,
\item a defuzzification method that provides an output value.
\end{itemize}

Our strategy is based on a 2-input/1-output inference scheme, with the fuzzy variables $k, \alpha$ and $w$. The output is the value of $w$. The first input is the current step $k$ of the algorithm. The second input is the percentage distance 
\[\alpha=\frac{(f(\textbf x_p(k)) - f(\mathbf{b}))\cdot 100}{f(\mathbf{b})}\]
between the objective function $f(\textbf x_p(k))$ evaluated on the current position of the particle $p$ at step $k$ and the objective function $f(\mathbf{b})$ evaluated on the global best position at step $k$. When $\alpha$ is small, it may be convenient to move the particle around its current position, corresponding to a low value of $w$. On the contrary, when $\alpha$ is large it may be sensible moving the particle far from its current position, corresponding to a high value of $w$.

During the execution of the PSO algorithm, the value of $w$ can be evaluated at each step $k$ in the following way. 

\begin{figure}[htp]
\centering
\subfloat[]
{
    \includegraphics[height=4cm,keepaspectratio]{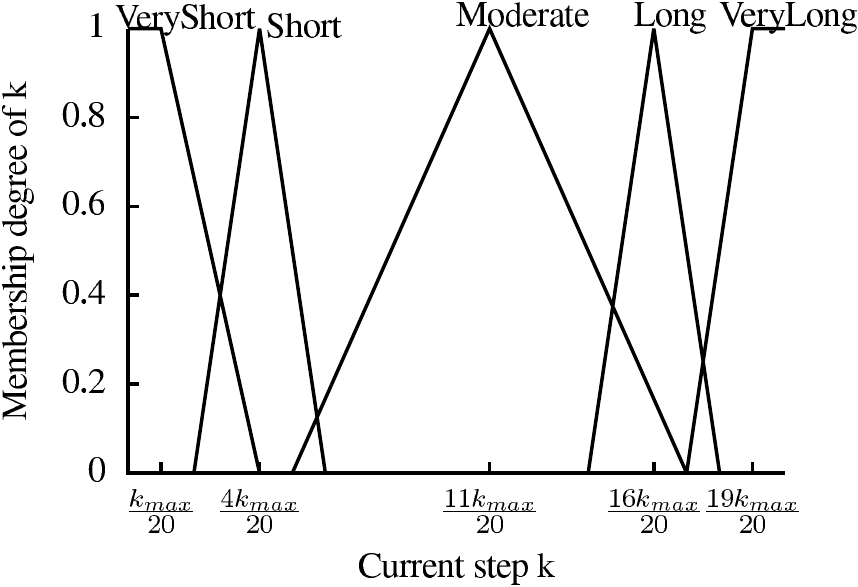}
    \label{fig:fuzzystep}
}%
\subfloat[]
{
    \includegraphics[height=4cm,keepaspectratio]{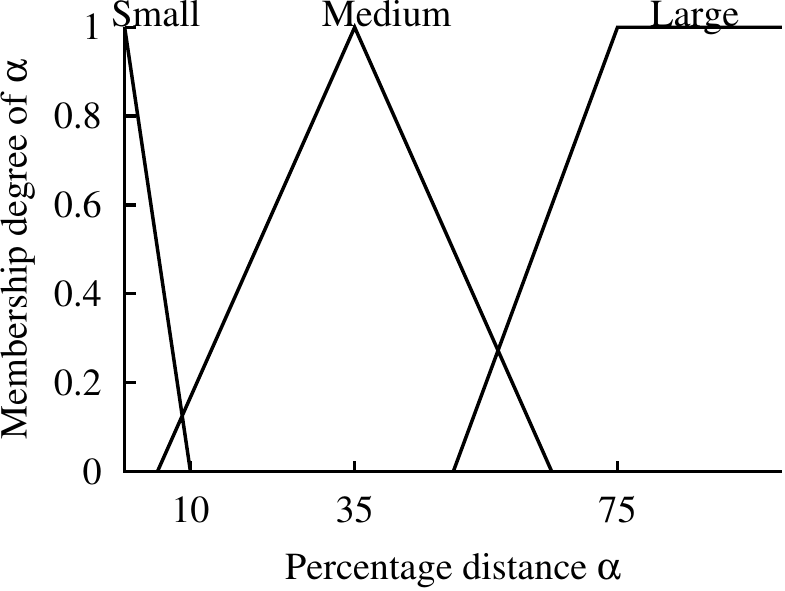}
    \label{fig:fuzzydist}
}\\
\subfloat[]
{
    \includegraphics[height=4cm,keepaspectratio]{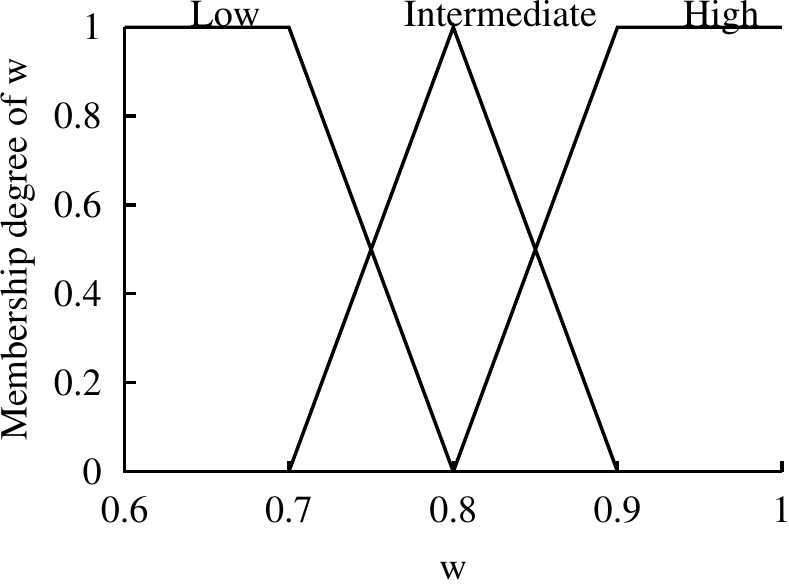}
    \label{fig:fuzzyw}
}
\caption{Membership degrees of the current step $k$, as a function of the value of $k$, for different fuzzy sets: VeryShort, ..., VeryLong~\protect\subref{fig:fuzzystep}, membership degrees of the percentage distance $\alpha$, as a function of the value of $\alpha$, for different fuzzy sets: Small, Medium, Large~\protect\subref{fig:fuzzydist} and membership degrees of $w$ for different fuzzy sets: Low, Intermediate, High~\protect\subref{fig:fuzzyw}.}
\label{fig:fuzzysdw}
\end{figure}

A membership degree is assigned to the current step $k$, as shown in Figure~\ref{fig:fuzzystep} for different fuzzy sets: VeryShort, Short, Moderate, Long, VeryLong. A fuzzy set characterizes the current step $k$ in terms of the membership degree of $k$, according to the following rules (the choice of values are expert based): 
\begin{enumerate} \small
\item if $k\leq\frac{4k_{max}}{20}$, then current step $k$ is VeryShort,
\item if $\frac{2k_{max}}{20}\leq k\leq\frac{6k_{max}}{20}$, then current step $k$ is Short,
\item if $\frac{5k_{max}}{20}\leq k\leq\frac{17k_{max}}{20}$, then current step $k$ is Moderate,
\item if $\frac{14k_{max}}{20}\leq k\leq\frac{18k_{max}}{20}$, then current step $k$ is Long,
\item if $k\geq\frac{17k_{max}}{20}$, then current step $k$ is VeryLong.
\end{enumerate}

For example, looking at Figure \ref{fig:fuzzystep} if the value of $k$ is less or equal than $\frac{k_{max}}{20}$, then $k$ belongs to the fuzzy set VeryShort with membership degree 1. If $\frac{2k_{max}}{20}\leq k \leq \frac{4k_{max}}{20}$, then $k$ belongs to the fuzzy sets VeryShort and Short with membership degrees given by the ordinates of the corresponding curves.

Similarly, a membership degree is assigned to $\alpha$, as shown in Figure~\ref{fig:fuzzydist}, according to the following rules:
\begin{enumerate} \small
\item if $\alpha\leq 10\%$, then distance is Small,
\item if $5\% \leq \alpha \leq 65\%$, then distance is Medium,
\item if $\alpha \geq 50\%$, then current step is Large,
\end{enumerate}
where Small, Medium and Large are the labels of the fuzzy sets that characterize the distance variable $\alpha$. Let us observe that if $\alpha<0$, then we set $\alpha=0$ as input. Indeed, in this case the current particle has improved the global best and consequently we would like to perform the search around this position giving to $\alpha$ the maximum membership degree 1 at the fuzzy set Small.

Finally, the membership degree of $w$ is obtained, as shown in Figure \ref{fig:fuzzyw}, according to the following rules (the choice of values are expert based): 
\begin{enumerate} \small
\item if $0.6 \leq w \leq 0.8$, then distance is Low,
\item if $0.7 \leq w \leq 0.9$, then distance is Intermediate,
\item if $0.8 \leq w \leq 1$, then current step is High,
\end{enumerate}
where Low, Intermediate and High are the labels of the fuzzy sets related to $w$. 
 
The derived membership functions (depicted by triangular or trapezoidal shapes) reflect expert--based choices. 
\begin{table}[!htb]\footnotesize
\centering
\caption{Inference system to derive the value of $w$.}
\begin{tabular}{|c|c||c|}
\hline 
$\bm{k}$ & $\bm{\alpha}$ & $\bm{w}$ \cr 
\hline \hline
 VeryShort & Small & Intermediate \cr
\hline
VeryShort &  Medium or Large & High \cr
\hline
Short & Small & Low \cr
\hline
Short & Medium or Large & High \cr
\hline
Moderate & Small & Low \cr
\hline
Moderate & Medium & Intermediate \cr
\hline
Moderate & Large & High \cr
\hline
Long & Small & Low \cr
\hline  
Long & Medium or Large & Intermediate \cr
\hline
VeryLong & Small or Medium & Low \cr
\hline
VeryLong & Large & Intermediate \cr
\hline
\end{tabular}
\label{table:fuzzyinfrules}
\end{table}
The inference system, used by the inference engine to derive the value of $w$, is based on the eleven rules shown in Table~\ref{table:fuzzyinfrules}.

The inference engine is the basic Mamdani model \cite{Mam} with if--then rules, minimax set--operations, sum for composition of activated rules, and defuzzification based on the centroid method that provides the abscissa of the barycentre of the fuzzy set composed according to the activated rules. 

Finally, the values of $c_1$ and $c_2$ can be derived as a function of $w$ from~\eqref{c1c2}.

\section{Variant 3: Bayesian approach to change the position of particles} \label{sec:bayes}

In this section, we propose a new strategy to change the position of the particles with respect to the PSO algorithm shown in Algorithm~\ref{alg:BasicPSOalgorithm}.
Specifically, in Subsection~\ref{sec:bayesianalgorithm} we derive the new algorithm based on the Kalman filter, and in Subsection~\ref{sec:parbayesianalgorithm} we describe the new parameters of the algorithm and we analyze how they impact on the velocity of the particles. 

\subsection{Bayesian PSO algorithm}\label{sec:bayesianalgorithm}

In \cite{Der2}, the authors highlighted that the Kalman filter can be taken back to the Bayes theorem, considering random variables with normal distributions.

Let $L|X$, $X$ and $X|L$ be random variables denoting, in terms of the Kalman filter, respectively, the current measurement (conditional on unknown parameter $X$), the prior estimation and the posterior estimation of the system's state. 
The Bayes theorem states
\[ f_{X|L}(x|l) = \frac{f_{L|X}(l|x)f_X(x)}{\int_{-\infty}^{+\infty}f_{L|X}(l|u)f_X(u)du}, \]
where $f_{X|L}$ is the posterior density, $f_{X}$ the prior density, $f_{L|X}$ the likelihood and the denominator is a normalization factor. 

Assume that $L|X$ and $X$ have the normal distributions $\mathcal N(\mu_{L|X},\sigma_{L|X}^2)$ and $\mathcal N(\mu_{X},\sigma_{X}^2)$, respectively. Then, it is well--known that the posterior estimation $X|L$ has the normal distribution $\mathcal N(\mu_{X|L},\sigma_{X|L}^2)$, where
\[ \mu_{X|L}=\frac{\sigma_{L|X}^2\mu_X+\sigma_X^2\mu_{L|X}}{\sigma_{L|X}^2+\sigma_X^2}, \quad \sigma_{X|L}^2=\left( \frac{1}{\sigma_{L|X}^2}+\frac{1}{\sigma_X^2} \right)^{-1}. \]

Now, we define the new position of the particles at step $k$, for $k \ge 1$, of the PSO algorithm, in terms of the Kalman filter. 

For each particle $p=1,\dots,D$ and for each dimension $i=1,\dots,d$, we consider that: 
\begin{itemize}
\item the global best position $b_i(k-1)$ (for the sake of brevity denoted as $b_i$)  is the current measurement $B|Y$, 
\item the position $x_{pi}(k-1)$ is the prior estimation $Y=X|L$, and 
\item the new position $x_{pi}(k)$ is the posterior estimation $Y|B$.
\end{itemize}
The prior estimation $Y$ is obtained by the posterior estimation $X|L$ of the new position $x_{pi}(k)$, considering that:
\begin{itemize}
\item the personal best position $l_{pi}(k-1)$ (for the sake of brevity denoted as $l_{pi}$) is the current measurement $L|X$, 
\item the position $x_{pi}(k-1)$ is the prior estimation $X$, and 
\item the new position $x_{pi}(k)$ is the posterior estimation $X|L$.
\end{itemize}
Thus, observing that $L|X \sim \mathcal N(l_{pi},\sigma_{L|X}^2)$ and $X \sim \mathcal N(x_{pi}(k-1),\sigma_{pi}^2(k-1))$, 
the posterior density of the random variable $Y=X|L$ is the normal distribution $\mathcal N(\mu_{X|L}(k),\sigma_{X|L}^2(k))$ derived by the Bayes theorem as
\begin{align*} 
  \mu_{X|L}(k)&=\frac{\sigma_{L|X}^2 x_{pi}(k-1)+\sigma_{pi}^2(k-1) l_{pi}}{\sigma_{L|X}^2+\sigma_{pi}^2(k-1)},\\
  \sigma_{X|L}^2(k)&=\left( \frac{1}{\sigma_{L|X}^2}+\frac{1}{\sigma_{pi}^2(k-1)} \right)^{-1}.
\end{align*}
Finally, observing that $B|Y \sim \mathcal N(b_i,\sigma_{B|Y}^2)$ and $Y \sim \mathcal N(\mu_{X|L}(k),\sigma_{X|L}^2(k))$,  
the posterior density of the random variable $Y|B$ is the normal distribution $\mathcal N(\mu_{Y|B}(k),\sigma_{Y|B}^2(k))$ derived by the Bayes theorem as
\begin{align} 
  \mu_{Y|B}(k)&=\frac{ x_{pi}(k-1)+ \delta_{L|X}(k-1) l_{pi} + \delta_{B|Y}(k-1) b_i}{1 + \delta_{L|X}(k-1) + \delta_{B|Y}(k-1)}, \label{eq:muYB}\\
  \sigma_{Y|B}^2(k)&= \frac{\sigma_{pi}^2(k-1)}{1 + \delta_{L|X}(k-1) + \delta_{B|Y}(k-1)}, \label{eq:sigma2YB}
\end{align}
where 
\begin{align}
\delta_{L|X}(k-1)=\frac{\sigma_{pi}^2(k-1)}{\sigma_{L|X}^2} \text{\quad and \quad} 
\delta_{B|Y}(k-1)=\frac{\sigma_{pi}^2(k-1)}{\sigma_{B|Y}^2}. \label{eq:deltaLXBY}
\end{align} 

Using this Bayesian approach to change the position of the particles, 
the value for $x_{pi}(k)$, at each step $k$, for each particle $p$ and each dimension $i$, is derived from the random variable $Y|B$ with normal distribution $\mathcal N(\mu_{Y|B}(k),\sigma_{Y|B}^2(k))$ as
\begin{align} 
  x_{pi}(k)&=Y|B. \label{eq:xpik}
\end{align}
The mean $\mu_{Y|B}(k)$ is derived from equation~\eqref{eq:muYB} replacing recursively $x_{pi}(k-1)$ by the value of $\mu_{Y|B}(k-1)$ (the mean of $Y|B$ at previous step $k-1$).
The variance $\sigma_{Y|B}^2(k)$ is derived from equation~\eqref{eq:sigma2YB}, replacing recursively $\sigma_{pi}^2(k-1)$ by the value of $\sigma_{Y|B}^2(k-1)$ (the variance of $Y|B$ at previous step $k-1$), when the position $\mathbf{x}_p(k-1)$ is the new global best, i.e., $f(\mathbf{x}_p(k-1))<f(\mathbf{b})$, otherwise, when $f(\mathbf{x}_p(k-1))\ge f(\mathbf{b})$, $\sigma_{Y|B}^2(k)$ is considered equal to $\sigma_{Y|B}^2(k-1)$.

In this way, considering $h(k)$ the number of steps with a new global best until the step $k$, for $k\ge 1$, i.e.,
\begin{align} 
  h(k)&=
\begin{cases}
  h(k-1)+1 & \text{if } f(\mathbf{x}_p(k-1))<f(\mathbf{b})\\
  h(k-1),              & \text{otherwise.}
\end{cases}
\end{align}
with $h(0)=0$, then the closed formulas of $\sigma_{Y|B}^2(k)$, $\delta_{L|X}(k)$ and $\delta_{B|Y}(k)$ for $k \ge 1$ are given by
\begin{align} 
\sigma_{Y|B}^2(k)&=\frac{\sigma_{pi}^2(0)}{1+h(k) (\delta_{L|X}(0)+\delta_{B|Y}(0))}=\frac{\sigma_{pi}^2(0)}{1 + h(k) \left(\frac{\sigma_{pi}^2(0)}{\sigma_{L|X}^2}+\frac{\sigma_{pi}^2(0)}{\sigma_{B|Y}^2}\right)}, \label{eq:sigma2YBk} \\
\delta_{L|X}(k)&=\frac{\delta_{L|X}(0)}{1+h(k) (\delta_{L|X}(0)+\delta_{B|Y}(0))}=\frac{\sigma_{pi}^2(0)}{\sigma_{L|X}^2+h(k) \sigma_{pi}^2(0) \left(1+\frac{\sigma_{L|X}^2}{\sigma_{B|Y}^2}\right)}, \label{eq:deltaLXk} \\
\delta_{B|Y}(k)&=\frac{\delta_{B|Y}(0)}{1+h(k) (\delta_{L|X}(0)+\delta_{B|Y}(0))}=\frac{\sigma_{pi}^2(0)}{\sigma_{B|Y}^2+h(k) \sigma_{pi}^2(0) \left(1+\frac{\sigma_{B|Y}^2}{\sigma_{L|X}^2}\right)}. \label{eq:deltaBYk}
\end{align} 
The function $h(k)$ is monotone non-decreasing and $h(k) \le k$.
 
The dimension $i$ of the velocity of the particle $p$ at step $k$, denoted by $v_{pi}(k)$, is derived as the difference of two normal random variables $x_{pi}(k)-x_{pi}(k-1)$ (as shown in Appendix~\ref{appendix}), that is a normal random variable with distribution $\mathcal N(\mu_{V}(k),\sigma_{V}^2(k))$~\cite{T02b}, where
\begin{align} 
  &\mu_{V}(k)=\mu_{Y|B}(k)-\mu_{Y|B}(k-1) \nonumber \\
  &=\frac{\delta_{L|X}(k-1)}{1 + \delta_{L|X}(k-1) + \delta_{B|Y}(k-1)} (l_{pi}(k-1) - x_{pi}(k-1)) \nonumber \\
  &\hphantom{= }+ \frac{\delta_{B|Y}(k-1)}{1 + \delta_{L|X}(k-1) + \delta_{B|Y}(k-1)} (b_i(k-1) - x_{pi}(k-1)), \label{eq:muYBv} \\
  &\sigma_{V}^2(k)=\sigma_{Y|B}^2(k)+\sigma_{Y|B}^2(k-1) \nonumber \\
  &\phantom{\sigma_{V}^2(k)}=\frac{(2+\delta_{L|X}(k-1) + \delta_{B|Y}(k-1))\sigma_{pi}^2(k-1)}{1+\delta_{L|X}(k-1) + \delta_{B|Y}(k-1)}. \label{eq:sigma2YBv}
\end{align}
In~\eqref{eq:muYBv} the velocity $v_{pi}(k)$ results from the sum of the cognitive and social components of the PSO. 
Thus, as soon as the current position of a particle converges to the current personal and to the current global best position, the mean of the velocity $\mu_{V}(k)$ becomes $0$ and the expected new position of the particle tends to be equal to the current one. 
In this case, although the mean of the new velocity is $0$, a new position, that can be also significantly different from the current position, can occur with a probability that depends on the value of the variance $\sigma_{Y|B}^2(k-1)$.

A different formula for $\mu_{V}(k)$ can be derived from~\eqref{eq:muYBv} (as described in Appendix~\ref{appendix}), showing that the velocity $v_{pi}(k)$ can result from the sum of three different components (the momentum, the difference between the local best positions of the last two steps and the difference between the global best positions of the last two steps), obtaining
\begin{align} 
  \mu_{V}(k)&=\frac{1 - \delta_{L|X}(k-1) - \delta_{B|Y}(k-1)}{1 + \delta_{L|X}(k-1) + \delta_{B|Y}(k-1)} \mu_{V}(k-1) \nonumber \\
  &\hphantom{= }+\frac{\delta_{L|X}(k-1)}{1 + \delta_{L|X}(k-1) + \delta_{B|Y}(k-1)} (l_{pi}(k-1)-l_{pi}(k-2)) \nonumber \\
  &\hphantom{= }+\frac{\delta_{B|Y}(k-1)}{1 + \delta_{L|X}(k-1) + \delta_{B|Y}(k-1)} (b_i(k-1)-b_i(k-2)). \label{eq:muYBv3c}
\end{align}
Equation~\eqref{eq:muYBv3c} shows that, when the mean $\mu_{V}(k)$ of the velocity of a particle is greater than $0$, the particle changes direction with respect to the current velocity (moment or inertia component) only when new current personal or global best positions occur. Otherwise, when personal or global best positions do not change then the magnitude of new velocity is obtained as a fraction of the current velocity.

When the $i$-th dimension of the research space is bounded by $x^{min}_i$ and $x^{max}_i$, the value $x_{pi}(k)$ is obtained from the random variable $Y|B$ with normal distribution $\mathcal N(\mu_{Y|B}(k),\sigma_{Y|B}^2(k))$ defined in~\eqref{eq:muYB} and~\eqref{eq:sigma2YB}, using the formula
\begin{align} 
  x_{pi}(k)&=x^{min}_i + (Y|B - \bar{x}^{min} ) \frac{ x^{max}_i-x^{min}_i}{ \bar{x}_i^{max}-\bar{x}_i^{min} }, \label{eq:xpikbounded}
\end{align} 
where
\begin{align} 
  \bar{x}_i^{min}=&\min\{x^{min}_i,\mu_{Y|B}(k)-3 \sigma_{Y|B}(k)\}, \label{eq:barxmin} \\
  \bar{x}_i^{max}=&\max\{x^{max}_i,\mu_{Y|B}(k)+3 \sigma_{Y|B}(k)\}. \label{eq:barxmax}
\end{align} 

In this way, values sampled from $Y|B$ that are in the interval $[\bar{x}^{min}_i,\bar{x}^{max}_i]$ are adjusted with values in the interval $[x^{min}_i,x^{max}_i]$ of the search variable.
The interval $[\bar{x}^{min}_i,\bar{x}^{max}_i]$ is defined in~\eqref{eq:barxmin} and~\eqref{eq:barxmax} such that the probability that $Y|B$ lies out of this interval is small. This probability is less than the probability that $Y|B$ is out of the interval $[\mu_{Y|B}(k)-3 \sigma_{Y|B}(k),\mu_{Y|B}(k)+3 \sigma_{Y|B}(k)]$, that is equal to $0.003$~\cite{T02b}.
The values sampled from $Y|B$ that are out of the interval $[\bar{x}^{min}_i,\bar{x}^{max}_i]$ are adjusted with the bounds $x^{min}_i$ or $x^{max}_i$ of the search variable.

Algorithm~\ref{alg:PSOalgorithmBayesianxpi} describes formally in pseudo–code the detailed steps to update the position of each particle at step $k$, with $k\ge 1$.
\begin{algorithm}[!htb]
  \DontPrintSemicolon
  \KwData{Step $k$, $\sigma_{B|Y}^2$, $\sigma_{L|X}^2$, $\mu_{Y|B}(k-1)$ and $\sigma_{Y|B}^2(k-1)$ for $p$ and $i$\; 
    $f(\mathbf{x}_p(k-1))$, $f(\mathbf{b}(k-1))$\;
    Bounds $x_i^{min}$ and $x_i^{max}$ of the $i$-th dimension of the particle's position\;
  }
  \KwResult{
    $x_{pi}(k)$, $\mu_{Y|B}(k)$ and $\sigma_{Y|B}^2(k)$ for $p$ and $i$
  }
  \Begin{
    \For(\tcp*[f]{for each particle}){$p=1,\dots,D$}{
      new\_global\_best=$f(\mathbf{x}_p(k-1)) < f(\mathbf{b}(k-1))$ \nllabel{algl:fxvsfb}\;
      \For(\tcp*[f]{for each search variable}){$i=1,\dots,d$}{
        $x_{pi}(k-1)\gets \mu_{Y|B}(k-1)$ \;
        $\sigma_{pi}^2(k-1)\gets \sigma_{Y|B}^2(k-1)$ \;
        $\mu_{Y|B}(k)\gets$ the right expression of equation~\eqref{eq:muYB} \;
        \If{new\_global\_best \nllabel{algl:isnewglobalbest}}
        {$\sigma_{Y|B}^2(k)\gets$ the right expression of equation~\eqref{eq:sigma2YB}}
        \Else{$\sigma_{Y|B}^2(k)\gets \sigma_{Y|B}^2(k-1)$}
        $x \gets \mathcal N(\mu_{Y|B}(k),\sigma_{Y|B}^2(k))$ \tcp*{new normal random value}
        $\bar{x}_i^{min} \gets \min\{x^{min}_i,\mu_{Y|B}(k)-3 \sigma_{Y|B}(k)\}$ \;
        $\bar{x}_i^{max} \gets \max\{x^{max}_i,\mu_{Y|B}(k)+3 \sigma_{Y|B}(k)\}$ \;
        \If{$x < \bar{x}_i^{min}$}{$x_{pi}(k)\gets x_i^{min} $}
        \ElseIf{$x > \bar{x}_i^{max}$}{$x_{pi}(k)\gets x_i^{max} $}
        \Else{
          $x_{pi}(k) \gets x^{min}_i + (x - \bar{x}^{min} ) \frac{ x^{max}_i-x^{min}_i}{ \bar{x}_i^{max}-\bar{x}_i^{min}}$
        }
      }
    }
}
\caption{Bayesian approach to update the position of each particle at step $k$, with $k \ge 1$}
\label{alg:PSOalgorithmBayesianxpi}
\end{algorithm}
The condition at line~\ref{algl:fxvsfb} verifies if the position of the particle $p$ at the previous step $k-1$ is the new global best position for the step $k$, i.e., if it is lower than the global best position obtained at the previous step $k-1$. This condition, that does not depend on the value of $i$, is used at line~\ref{algl:isnewglobalbest} to decide the current value of the variance $\sigma_{Y|B}^2(k)$, for each $i$. 

\subsection{Parameters and analysis of the Bayesian PSO algorithm}\label{sec:parbayesianalgorithm}

As shown in~\eqref{eq:muYB}, the new position of a particle is determined by means of a weighted mean among the previous position of the particle, the current personal best position and the current global best position. The weights of each component are determined by $\delta_{L|X}(k-1)$ and $\delta_{B|Y}(k-1)$, which are defined as a function of the variances $\sigma_{B|Y}^2$, $\sigma_{L|X}^2$ and $\sigma_{pi}^2(0)$. These variances represent the uncertainties associated to the current measurement at each step of $b_i$ (conditioned to the prior estimation of $l_i$), the current measurement at each step of $l_i$ (conditioned to the prior estimation of the position $x_{pi}$) and the posterior estimation at step $1$ of the new position $x_{pi}(1)$, respectively.
Parameters $\sigma_{B|Y}^2$ and $\sigma_{L|X}^2$ correspond to the acceleration coefficients of the basic PSO algorithm, since they control the impact of the cognitive and social components on the velocity of the particle.
Parameter $\sigma_{pi}^2(0)$ adjusts the position of a particle with respect to the new value of the mean $\mu_{Y|B}(k)$.
From~\eqref{eq:deltaLXBY} and~\eqref{eq:muYBv}  it follows that, the ratio of the weights of the cognitive and social components is
\begin{align}
\frac{\sigma_{B|Y}^2}{\sigma_{L|X}^2}. \label{eq:ratioWeightCognitiveSocial}
\end{align}
Thus, when $\sigma_{L|X}^2 < \sigma_{B|X}^2$, the impact of the cognitive component on the velocity of the particle is greater than of the impact of the social component, like in the basic PSO algorithm when $c_1>c_2$. In this case, the particle is more attracted to its own personal best position, resulting in excessive wandering~\cite{E07}. On the other hand, if $\sigma_{L|X}^2 > \sigma_{B|X}^2$ then, like in the  basic PSO algorithm when $c_1<c_2$, the particle is strongly attracted to the global best position, causing particles to converge prematurely towards global best position~\cite{E07}.
Considering different values of $\sigma_{L|X}^2$ and $\sigma_{B|Y}^2$ while maintaining the same ratio, different values for $\mu_{V}(k)$ can be obtained without changing the relative impacts of the cognitive and social components.

From~\eqref{eq:muYBv} it follows that, for $k \ge 1$:
\begin{itemize}
\item the mean of the new position $\mu_{Y|B}(k)$ remains equal to the mean of the current position $\mu_{Y|B}(k-1)$, i.e., the mean of the velocity $\mu_{V}(k)$ tends to $0$, as $\delta_{L|X}(k-1)$ and $\delta_{B|Y}(k-1)$ approach to $0$, and 
\item $\mu_{Y|B}(k)$ tends to $l_{pi}(k-1)$ or $b_i(k-1)$ when, respectively, $\delta_{L|X}(k-1)$ or $\delta_{B|Y}(k-1)$ tend to infinity.
\end{itemize}
In fact, $\delta_{L|X}(k)$ or $\delta_{B|Y}(k)$ can assume very high values only when $k=0$.
Definitions~\eqref{eq:deltaLXBY} for $k=0$ imply that $\delta_{L|X}(0)$ or $\delta_{B|Y}(0)$ vary from $0$ to infinity as $\sigma_{L|X}^2$ or $\sigma_{B|Y}^2$, respectively, vary from infinity to 0, or as $\sigma_{Y|B}^2(0)$ varies from $0$ to infinity.

From the closed formulas~\eqref{eq:deltaLXk} and~\eqref{eq:deltaBYk} we derive that, depending on the initial value assigned to $\sigma_{pi}^2(0)$, for $h(k) \ge 1$
\begin{align} 
0 < \delta_{L|X}^2(k) < \frac{1}{ h(k) \left(1+\frac{\sigma_{L|X}^2}{\sigma_{B|Y}^2} \right)} \text{, \quad}
0 < \delta_{B|Y}^2(k) < \frac{1}{ h(k) \left(1+\frac{\sigma_{B|Y}^2}{\sigma_{L|X}^2} \right)} \label{eq:deltaLXBYkconstr}
\end{align}
where $\delta_{L|X}(k)$ and $\delta_{B|Y}(k)$ become $0$ when $\sigma_{pi}^2(0)$ tends to $0$, and they tend to their respective (less than $1$) upper bounds when $\sigma_{pi}^2(0)$ approaches infinity. 
Thus, the maximum value that $\delta_{L|X}^2(k)$ and $\delta_{B|Y}^2(k)$ can reach for $h(k) \ge 1$ depends on the ratio between the values of $\sigma_{L|X}^2$ and $\sigma_{B|Y}^2$ and on the value of $h(k)$, independently from the value of $\sigma_{pi}^2(0)$.

Replacing the bounds of equation~\eqref{eq:deltaLXBYkconstr} in equation~\eqref{eq:muYBv} we get
\begin{align} 
  &0 \le \mu_{V}(k) \le \frac{1}{h(k)\left(1+\frac{\sigma_{L|X}^2}{\sigma_{B|Y}^2}\right)} (l_{pi}(k-1) - x_{pi}(k-1)) \nonumber \\
  &\hphantom{0 \le \mu_{V}(k) \le }+ \frac{1}{h(k)\left(1+\frac{\sigma_{B|Y}^2}{\sigma_{L|X}^2}\right)} (b_i(k-1) - x_{pi}(k-1)), \label{eq:muYBvBound}
\end{align}
where, for $h(k) \ge 1$, the mean $\mu_{V}(k)$ of the velocity becomes $0$, when $\sigma_{pi}^2(0)$ approaches $0$, and it tends to the upper bound of~\eqref{eq:muYBvBound}, when $\sigma_{pi}^2(0)$ approaches infinity.

From~\eqref{eq:deltaLXk} and~\eqref{eq:deltaBYk} it follows also that, for $h(k) \ge1$, $\delta_{L|X}^2(k)$ and $\delta_{B|Y}^2(k)$ have opposite trends, with $\delta_{L|X}^2(k)$ decreasing and $\delta_{B|Y}^2(k)$ increasing for increasing values of $\sigma_{L|X}^2$, and with $\delta_{L|X}^2(k)$ increasing and $\delta_{B|Y}^2(k)$ decreasing for increasing values of $\sigma_{B|Y}^2$. This implies that, when the two (cognitive and social) components of the velocity in~\eqref{eq:muYBv} are positive and for a fixed value of $\sigma_{pi}^2(0)$ and $h(k)$, the mean $\mu_{V}(k)$ of the velocity is greater than a positive lower bound; it is not dependent on the values of $\sigma_{L|X}^2$ and $\sigma_{B|Y}^2$, but only on the value of $\sigma_{pi}^2(0)$ and $h(k)$.

Thus, when the cognitive and social components of the velocity in~\eqref{eq:muYBv} are positive, $\mu_{V}(k)$  for $h(k) \ge1$ is limited by a positive lower bound, that depends only on $h(k)$ and $\sigma_{pi}^2(0)$, and by the upper bound of~\eqref{eq:muYBvBound}, that depends only on $h(k)$ and on the ratio between $\sigma_{L|X}^2$ and $\sigma_{B|Y}^2$.

Depending on the initial value assigned to $\sigma_{pi}^2(0)$, from~\eqref{eq:sigma2YBk} we have for $h(k) \ge1$
\begin{align} 
0 < \sigma_{Y|B}^2(k) < \frac{1}{ h(k) \left(\frac{1}{\sigma_{L|X}^2}+\frac{1}{\sigma_{B|Y}^2} \right)} = 
\frac{\sigma_{L|X}^2}{ h(k) \left(1+\frac{\sigma_{L|X}^2}{\sigma_{B|Y}^2} \right)} \label{eq:sigma2YBkconstr}
\end{align} 
where $\sigma_{Y|B}^2(k)$ becomes $0$, when $\sigma_{pi}^2(0)$ approaches $0$, and it tends to the upper bound of~\eqref{eq:sigma2YBkconstr}, when $\sigma_{pi}^2(0)$ approaches infinity. 
Thus, the maximum value that the variance $\sigma_{Y|B}^2(k)$ and $\sigma_V^2(k)$ can reach is limited by the values of the parameters $\sigma_{L|X}^2$ and $\sigma_{B|Y}^2$, independently from the value of $\sigma_{pi}^2(0)$.
From~\eqref{eq:sigma2YBk} we derive that, for values of $\sigma_{L|X}^2>>\sigma_{pi}^2(0)$ or $\sigma_{B|Y}^2>>\sigma_{pi}^2(0)$ and for a limited value of $h(k)$, the variance $\sigma_{Y|B}^2(k)$ can be approximated by the constant $\sigma_{pi}^2(0)$.

As an example of parameter setting, considering $\sigma_{L|X}^2=\sigma_{B|Y}^2$ and $\sigma_{pi}^2(0)>>\sigma_{L|X}^2$, we get for $h(k) \ge1$
\begin{align*}
  \delta_{L|X}^2(k)=\delta_{B|Y}^2(k)\simeq \frac{1}{2h(k)} \text{\quad and \quad}
  \sigma_{Y|B}^2(k)\simeq \frac{\sigma_{L|X}^2}{2h(k)},
\end{align*}
where the right side expressions of the two equations are derived from the upper bounds of~\eqref{eq:deltaLXBYkconstr} and~\eqref{eq:sigma2YBkconstr}, respectively.

When the number of new global best positions $h(k)$ of the algorithm increases, the uncertainties relative to the current personal best position and the current global best position decrease, since improved positions for the personal and global best should be reached. 
From the previous discussion about the trends of $\mu_{Y|B}(k)$ and from equations~\eqref{eq:sigma2YBk},~\eqref{eq:deltaLXk} and~\eqref{eq:deltaBYk} it follows that, for increasing values of $h(k)$, the mean $\mu_{Y|B}(k)$ tends to be constant and the variance $\sigma_{Y|B}^2(k)$ decreases to $0$.
In order to avoid that the mean $\mu_V(k)$ and the variance $\sigma_V^2(k)$ of the velocity decrease below a certain minimum value before the optimal position is reached, the following condition must be verified:
\begin{align}
  \sigma_{Y|B}^2(k_{max}) \ge \epsilon, \label{eq:sigma2YBepsilon}
\end{align}
where $\epsilon >0$ and $k_{max}$ is the maximum number of steps for which the optimal position is reached.
Given the values for $\sigma_{L|X}^2$, $\sigma_{B|Y}^2$, $k_{max}$ and $\epsilon$, then  
the minimum value for $\sigma_{pi}^2(0)$ for which the variance $\sigma_{Y|B}^2(k_{max}) \ge \epsilon$ is derived from~\eqref{eq:sigma2YBk} as
\begin{align} 
  \frac{\epsilon}{1-\epsilon h(k_{max}) \left(\frac{1}{\sigma_{L|X}^2}+\frac{1}{\sigma_{B|Y}^2} \right)}, \label{eq:sigma2YBkmin}
\end{align} 
with
\begin{align*} 
  \frac{1}{\sigma_{L|X}^2}+\frac{1}{\sigma_{B|Y}^2}<\frac{1}{\epsilon h(k_{max})}. %\label{eq:sigma2YBkgtzero}
\end{align*} 

For large values of $\sigma_{L|X}^2$ or $\sigma_{B|Y}^2$, formula~\eqref{eq:sigma2YBkmin} is approximated by $\epsilon$.
On the contrary, low values of $\sigma_{L|X}^2$ and $\sigma_{B|Y}^2$ imply smaller values for $\epsilon$ or $h(k_{max})$ and values for $\sigma_{pi}^2(0)$ higher than $\epsilon$.

\section{Evaluation of the proposed solutions} \label{sec:simulations}

In this section, we carry on an evaluation study to assess performance and quality indicators of our solutions via simulation. To this purpose, we set up two variants of the basic PSO algorithm by assembling the previously described features. Specifically, the two PSO variants we concentrate on are generated respectively by i) employing the novel Algorithm \ref{alg:PSOinitial} for the sample of the initial population and the fuzzy approach described in Section~\ref{sec:fuzzy} (referred in the following as PSOF variant), and ii) the novel Algorithm \ref{alg:PSOinitial} and the Algorithm \ref{alg:PSOalgorithmBayesianxpi} to change the position of particles (referred in the following as PSOB variant). 
Two kinds of evaluation are performed. First, the PSOF and PSOB solutions are compared with the basic PSO algorithm (indicated as PSOC) in terms of performance and quality indicators. Then, the new algorithms are employed in a simple but realistic use case in the electrical power sector taken from the literature and compared with a few other optimization solutions already adopted in the referred study. These analyses are detailed in the next two subsections.

\subsection{Analysis and comparison with the basic PSO}\label{sec:basic_only}

In this study, we consider the Rosenbrock function \cite{Rose} and the Griewank function \cite{Gri}, which are functions widely used to test the quality of optimization algorithms. In particular, we will use their multi--dimensional generalizations, i.e.,
\[ f(x_1,...,x_d)=\sum_{i=1}^{d-1}((1-x_i)^2+100(x_{i+1}-x_i^2)^2) \]
and
\[g(x_1,...,x_d)=1+\cfrac{1}{4000}\sum_{i=1}^dx_i^2-\prod_{i=1}^d\cos \cfrac{x_i}{\sqrt{i}},\]
respectively. Moreover, we consider -10 and 10 as the lower and upper bounds for each variable of the Rosenbrock function and -20, 20 as the lower and upper bounds for each variable of the Grienwank function.

The two indicators we analyzed are \emph{accuracy} $A$, which represents the ability of the algorithm to better approximate (or possibly reach) the optimal solution $\mathbf{b}^*$ within a given number of steps, and \emph{efficiency} $K$, which represents the promptness (in number of steps) of the algorithm in reaching the optimal solution $\mathbf{b}^*$ or a certain approximated value $\overline{a}$ of the objective function. $A$ and $K$ are random variables defined as
\begin{align*}
  A&=f(\mathbf{b}({k_{max}}))-f(\mathbf{b}^*), \\
  K&=\min \{k | f(\mathbf{b}(k)) \le \overline{a} \text{ and k is not limited}\}. 
\end{align*}
When the minimum of the objective function is $0$, i.e., $f(\mathbf{b}^*)=0$, then $A$ is equal to the minimum value of the objective function obtained by the PSO algorithm, i.e., $f(\mathbf{b}({k_{max}}))$.

To enrich the analyses and the comparison between the three variants, for both indicators we evaluated the minimum and maximum values obtained out of the set of simulations performed, as well as the mean value, i.e.,
\begin{align*}
  A^{min}_X&=\min\{A\}, A^{max}_X=\max\{A\} \quad\text{and}\quad A^{mean}_X=E[A],\\
  K^{min}_X&=\min\{K\}, K^{max}_X=\max\{K\} \quad\text{and}\quad K^{mean}_X=E[K],
\end{align*}
where $X$ is the name of the variant for wich the measure is derived: PSOC, PSOF or PSOB.

The maximum number of steps  $k_{max}$ and the swarm dimension $D$ are the two parameters that have been varied in the analyses, as specified in the tables summarizing the obtained results. Default values for the algorithms parameters, assumed in the simulations when not otherwise specified, are respectively: $k_{max}=150$ and $D=35$.

Moreover, in the simulations we use the following standard values for the PSOC parameters:
\[w=0.72984, \quad c_1=c_2=1.496172.\]
For the PSOF algorithm, the coefficients $w, c_1, c_2$ are determined at each step $k$ by using the fuzzy inference scheme described in Section \ref{sec:fuzzy}.

Finally, in the PSOB algorithm, the parameters $\sigma^2_{pi}(0)$, $\sigma^2_{L|X}$, $\sigma^2_{B|Y}$ are set by expert-based choice as follows:
\[\sigma^2_{p}(0)=\sigma^2_{L|X}=\cfrac{\lvert x_i^{max} - x_i^{min} \rvert}{2D},\quad \sigma^2_{B|Y}=\cfrac{\lvert x_i^{max} - x_i^{min} \rvert}{D},\]
for each component $i$ of the $p$-th particle.

\subsubsection{Accuracy evaluation results} \label{subsubsec:accuracy}  

Results of 100 simulation runs to evaluate the accuracy parameter are summarized in Table \ref{table:rose} for the PSOC, PSOF, PSOB algorithms used to minimize the Rosenbrock function with $d=3$. The Rosenbrock function with $d=3$ takes the minimum value, that is zero, at the point $(0,0,0)$.

\begin{table}[!htb]\footnotesize
\caption{Minimization of the Rosenbrock funtion with $d=3$ using PSOC, PSOF, PSOB algorithms, for 100 simulations.}
\centering 
\tabcolsep=0.15cm
\scalebox{0.82}{
\begin{tabular}{|c||c|c|c||c|c|c||c|c|c|}
\hline
& $\mathbf{A^{min}_{PSOC}}$ & $\mathbf{A^{mean}_{PSOC}}$ & $\mathbf{A^{max}_{PSOC}}$ & $\mathbf{A^{min}_{PSOF}}$ & $\mathbf{A^{mean}_{PSOF}}$ & $\mathbf{A^{max}_{PSOF}}$ & $\mathbf{A^{min}_{PSOB}}$ & $\mathbf{A^{mean}_{PSOB}}$ & $\mathbf{A^{max}_{PSOB}}$ \cr \hline \hline
\textbf{Default} & $4.8\cdot10^{-5}$ & 0.5730 & 8.7900 & 4.2$\cdot 10^{-7}$ & 0.04701 & 7.3425 & 0.0008 & 0.0339 &  0.1523 \cr \hline
$\mathbf{D=20}$ & 8.8$\cdot 10^{-5}$ & 0.9532 & 8.0209 & 8.0$\cdot 10^{-5}$ & 0.3920 & 9.2987  & 0.0033 & 0.1473 & 0.6662 \cr \hline
$\mathbf{D=50}$ & 1.1$\cdot 10^{-5}$ & 0.1951 & 4.9227 & 2.1$\cdot 10^{-7}$ & 0.0134 & 0.0645 & 2.1$\cdot 10^{-5}$ & 0.0127 & 0.1057  \cr \hline
$\mathbf{k_{max}=100}$ & $6.2632\cdot10^{-5}$ & 0.6534 & 5.8641 & 0.0002 & 0.0727 & 0.2868 & 0.0005 & 0.0552 & 0.2147 \cr \hline
$\mathbf{k_{max}=200}$ & 6.1$\cdot 10^{-6}$ & 0.3933 & 5.2793 & 4.3$\cdot 10^{-9}$ & 0.0284 & 0.1110 & 6.5$\cdot 10^{-5}$ & 0.0258 & 0.1179 \cr\hline
\end{tabular}
}
\label{table:rose}
\end{table}

We can observe that for the default values of the parameters (first row of the table), PSOF has a better accuracy than PSOB and PSOC with respect to the minimum value. %Indeed the minimum value of the Rosenbrock function is 4.2$\cdot 10^{-7}$ for the PSOF, 0.0008 for the PSOB, and $4.8\cdot10^{-5}$  for the PSOC. 
Instead, PSOB  shows a better accuracy than PSOF and PSOC for both $A^{mean}_{PSOB}$ and $A^{max}_{PSOB}$. 
Similar results are obtained for different values of the parameters $D$ and $k_{max}$. Although there is no definitive rank among the three variants of the PSO under analysis, the accuracy of the basic PSO version is never better than both the other two. Especially when focusing on $A^{mean}_{PSOC}$, it is constantly (and significantly) higher than $A^{mean}_{PSOF}$ and $A^{mean}_{PSOB}$, meaning lower accuracy.

Similarly, results of 100 simulation runs are summarized in Table \ref{table:gri} for the PSOC, PSOF, PSOB algorithms used to minimize the Griewank function with $d=5$, again to evaluate the accuracy indicator. The Griewank function with $d=5$ takes the minimum value, that is zero, in several different points.

\begin{table}[!htb]\footnotesize
\caption{Minimization of the Griewank funtion with $d=5$ using PSOC, PSOF, PSOB algorithms, for 100 simulations.}
\centering
\tabcolsep=0.15cm
\scalebox{0.82}{
\begin{tabular}{|c||c|c|c||c|c|c||c|c|c|}
\hline 
& $\mathbf{A^{min}_{PSOC}}$ & $\mathbf{A^{mean}_{PSOC}}$ & $\mathbf{A^{max}_{PSOC}}$ & $\mathbf{A^{min}_{PSOF}}$ & $\mathbf{A^{mean}_{PSOF}}$ & $\mathbf{A^{max}_{PSOF}}$ & $\mathbf{A^{min}_{PSOB}}$ & $\mathbf{A^{mean}_{PSOB}}$ & $\mathbf{A^{max}_{PSOB}}$ \cr \hline \hline
\textbf{Default} & 0.0004 & 0.0757 & 0.2234 & 2.2$\cdot10^{-6}$ & 0.0118 & 0.0236 & 0.0038 & 0.0071 &  0.0481 \cr \hline
$\mathbf{D=20}$ & 0.0072 & 0.0901 & 0.2280 & 0.0001 & 0.0839 & 0.2228  & 0.0005 & 0.0369 & 0.1073 \cr \hline
$\mathbf{D=50}$ & 9.1$\cdot 10^{-5}$ & 0.0626 & 0.2046 & 6.4$\cdot 10^{-10}$ & 8.2$\cdot10^{-5}$ & 0.0005 & 6.9$\cdot 10^{-5}$ & 0.0008 & 0.0035  \cr \hline
$\mathbf{k_{max}=100}$ & 0.0114 & 0.0941 & 0.2612 & 0.0007 & 0.0832 & 0.2490 & 0.0006 & 0.0075 & 0.0481 \cr \hline
$\mathbf{k_{max}=200}$ & 0.0003 & 0.0699 & 0.1626 & 4.7$\cdot 10^{-8}$ & 0.0007 & 0.0153 & 3.0$\cdot 10^{-5}$ & 0.0004 & 0.0495 \cr\hline
\end{tabular}
}
\label{table:gri}
\end{table}

Comments similar to those already made with reference to the previous table are applicable also to Table \ref{table:gri}. PSOC accuracy is always outperformed by either PSOB or PSOF (actually, by both, unless for the minimum value in the first row when the default parameters are used, for which PSOC is better than PSOB). Therefore, also in this use case the proposed variants PSOF and PSOB should be preferred.

\subsubsection{Efficiency evaluation results} \label{subsubsec:efficiency}

To assess the efficiency indicator, we considered the Griewank function only and calculated the (minimum, mean and maximum) number of steps required by each of the three PSO variants to reach a predefined value $\overline{a}$ for this function. Again, 100 simulations have been performed; each run is completed either when the specified value is obtained by the execution of the PSO algorithm, or a maximum of 150 steps is exceeded.  
Table \ref{table:efficiency} summarizes the simulation results. The values chosen for $\overline{a}$ are the minimum, mean and maximum values obtained in Table \ref{table:gri} by PSOC, in correspondence of the row ``Default'', that is, $A^{min}_{PSOC}=0.0004$, $A^{mean}_{PSOC}=0.0757$ and $A^{max}_{PSOC}=0.2234$, respectively. The values $0$ in the table indicate that the initial position randomly selected by the corresponding algorithms already provide an outcome better than $\overline{a}$. Instead, when ``not found'' appears, it means that $\overline{a}$ is not reached within 150 steps (maximum number of steps set in the simulations).
Looking at the results in Table~\ref{table:efficiency}, we can immediately appreciate the significantly better efficiency shown in general by the variants PSOF and PSOB. In fact, for both $\overline{a}=A^{mean}_{PSOC}$ and $\overline{a}=A^{max}_{PSOC}$, the number of steps requested by PSOF and PSOB is much lower than that required by PSOC, being indeed 0 for $\overline{a}=A^{max}_{PSOC}$. The only exception is PSOB for $\overline{a}=A^{min}_{PSOC}$, which is not surprising, since looking at Table \ref{table:gri}, PSOB shows lower accuracy than the minimum value obtained by PSOC, when default parameters are used (which is 0.0004). 

\begin{table}[!htb]\footnotesize
\caption{Minimization of the Griewank funtion with $d=5$ using PSOC, PSOF, PSOB algorithms, for 100 simulations.}
\centering 
\tabcolsep=0.15cm
\scalebox{0.82}{
\begin{tabular}{|c||c|c|c|c|c|c|c|c|c|}
\hline
& $\mathbf{K^{min}_{PSOC}}$ & $\mathbf{K^{mean}_{PSOC}}$ & $\mathbf{K^{max}_{PSOC}}$ & $\mathbf{K^{min}_{PSOF}}$ & $\mathbf{K^{mean}_{PSOF}}$ & $\mathbf{K^{max}_{PSOF}}$ & $\mathbf{K^{min}_{PSOB}}$ & $\mathbf{K^{mean}_{PSOB}}$ & $\mathbf{K^{max}_{PSOB}}$ \cr \hline \hline
$\overline{a}=\mathbf{A^{min}_{PSOC}}$  & 127 & 131 & not found & 16 & 28 & 50 & not found & not found & not found  \cr \hline
$\overline{a}=\mathbf{A^{mean}_{PSOC}}$ & 17 & 66 & 144 &  4  & 14 & 21 & 8 & 12 & 31 \cr \hline
$\overline{a}=\mathbf{A^{max}_{PSOC}}$  & 1 & 26 & 122 & 0 & 0 & 0 & 0 &  0  & 0  \cr \hline
\end{tabular}
}
\label{table:efficiency}
\end{table}

\subsection{Evaluation in a practical case study: optimization of electrical grids} \label{sec:electric}

\begin{figure}[!htb] 
\centering
\includegraphics[scale=2]{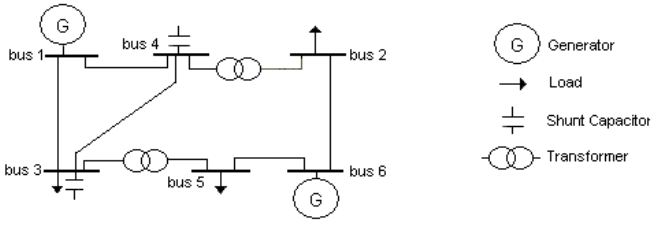}
\caption{IEEE--6 Bus System}
\label{fig:ieee6}
\end{figure}

In this section, we have evaluated our modified PSO algorithms when employed in  electrical grid optimizations, in terms of minimizing the total power loss over the lines. The evaluation is performed in terms of the accuracy indicator, as  introduced in \ref{sec:basic_only}. In particular, we focused on the IEEE--6 bus system in Figure \ref{fig:ieee6}, typically used as a simple case study to deal with optimization in electrical grids. A survey on the general optimal power flow problem (i.e., the problem to optimize electrical grids) can be found in \cite{Stevenson} and \cite{Stott}.
In \cite{Sin}, the authors have reported the results of the optimization of the IEEE--6 bus system by using different optimization algorithms: Genetic Algorithm (GA), Particle Swarm Optimization (PSO), Artificial Bee Colony (ABA), Differential Evolution (DE). The optimization has been implemented through minimizing the total power loss over the lines of the electrical grid, expressed as follows:
\[P_{loss}=\frac{1}{2}\sum_{i=1}^N\sum_{j=1}^n g_{i,j}(V_i^2+V_j^2-2V_iV_j\cos(\delta_i-\delta_j))\] 
where $N$ is the total number of bus (nodes), $V_i$ is the voltage magnitude at bus $i$, $\delta_i$ is the voltage angle at bus $i$, $g_{i,j}$ is the conductance of the line connecting busses $i$ and $j$. For optimization purposes, seven control variables have been used in these papers: $P_6$ active power of bus 6, $V_1$ and $V_6$ voltages of busses 1 and 6, $Q_3$ and $Q_4$ reactive power of busses 3 and 4, $T_{24}$ and $T_{35}$ tap changers of transformers over lines connecting busses 2,4 and 3,5. We have instantiated our PSOF and PSOB solutions to this problem, adopting the same control variables, which correspond to the search variables of our algorithms.  For the sake of comparison with the previously analysed optimization algorithms, the same range of values as in \cite{Sin} have been assigned to the (both control and state) variables of our algorithms; for brevity, they are omitted here.

In Table \ref{table:opf}, we report the results obtained in \cite{Sin} for the GA, PSO, ABA and DE algorithms, with the addition of the results that we have obtained using our modified algorithms PSOB and PSOF. These results (expressed in terms of min, mean and maximum total power loss) are obtained with $k_{max}=200$, $D=30$, running 30 simulation experiments. From this table, it can be observed that both PSOB and PSOF show the best minimum value for the total power loss.  Moreover, PSOB is the best also in terms of the average and maximum values, so ranking as the most accurate power flow optimization algorithm among the 6 considered alternatives for this scenario. Although dependent on the specific setting adopted in the set up of this case study, the observed results are encouraging to investigate more deeply the suitability of our proposed solutions to support power flow optimizations, which we plan as future work. 

\begin{table}[!htb]\footnotesize
\caption{OPF solutions obtained by using GA, PSO ABC, DE, PSOB, and PSOF for the 6--bus test system}
\centering 
\tabcolsep=0.15cm
\scalebox{0.82}{
\begin{tabular}{|c|c|c|c|}
\hline
\textbf{Method} & $\mathbf{P^{Min}_{loss}}$ & $\mathbf{P^{Mean}_{loss}}$ & $\mathbf{P^{Max}_{loss}}$ \cr \hline \hline
GA & 6.7747 & 6.9705 & 7.5292  \cr \hline
PSO & 6.7486 & 6.8425 & 7.1517 \cr \hline
ABA & 6.7361 & 6.7361 & 6.7364 \cr \hline
DE & 6.7361 & 6.7361 & 6.7368 \cr \hline
PSOF & 6.7329 & 6.7557 & 6.9042 \cr \hline
PSOB & 6.7329 & 6.7331 & 6.7333 \cr \hline
\end{tabular}
}
\label{table:opf}
\end{table}
 
For completeness on the presented case study, in Table \ref{table:optimal-cv} we report the optimal configuration for the control variables that yields the minimum value of the total power loss obtained by the PSOB and the PSOF algorithms.

\begin{table}[!htb]\footnotesize
\caption{Values of control variables corresponding to minimum solution obtained by PSOB and PSOF algorithms}
\centering 
\tabcolsep=0.15cm
\scalebox{0.82}{
\begin{tabular}{|c|c|}
\hline
\textbf{Control variable} & \textbf{Value}  \cr \hline \hline
$V_1$ & 1.1  \cr \hline
$V_6$ & 1.1  \cr \hline
$P_6$ & 27.6  \cr \hline
$Q_3$ & 43   \cr \hline
$Q_4$ & 27  \cr \hline
$T_{24}$ & 1.0475  \cr \hline
$T_{35}$ & 0.9975  \cr \hline
\end{tabular}
}
\label{table:optimal-cv}
\end{table}

\section{Conclusions} \label{sec:conclusions}

In this paper, the problem of performing non--linear optimization has been treated by means of the PSO algorithm. Moving from the original formulation of the PSO algorithm, improvements have been proposed on the sampling of the initial population, the setting of the algorithm's parameters and the method for moving the control variables in the search space.  Well consolidated techniques have been adopted to cover the critical points of the algorithm we worked on, namely the fuzzy logics and the Bayesian theory. A simulation study has been carried on to show the benefits of the new proposed PSO solutions with respect to the classical PSO formulation, using well--known testing functions like the Rosenbrock function and the Griewank function.  Moreover, a case study in the electrical field has been worked out, to show results on the suitability of our proposed solutions to support optimization needs in this context, in comparisons with other already adopted alternatives. The obtained results are encouraging, and we expect that the new features we introduced into the PSO algorithm are actually relevant in a variety of application contexts, especially those that are sensitive to the input sampling and the coverage of the search space. 

%This research effort opens the way to new or improved applications of the PSO--based algorithms to solve optimization problems in different engineering contexts, and possibly to be further enhanced by tackling other critical aspects of the original formulation. 

Several research extensions are foreseen. 
%A first direction would be to investigate further the benefits of the proposed variants by considering additional combinations of them, such as combining the solution for the sample of the initial population with the Bayesian technique to change the position of the particles.  
Further simulations devoted to improve the understanding of the sensitivity of the three compared solutions to the algorithms parameters would be undoubtedly a valuable direction to explore. 
Other refinements would be also interesting, such as focus on the setting of variances and introducing correlation factors in the Bayesian approach, as well as setting of expert--based choices in the fuzzy logics strategy. Of course, addressing more deep investigations on the usage of the proposed solutions in specific application contexts, such as the electrical power system already tackled in this paper, is another planned research line. Finally, practical support to the the selection of the most suited algorithm to solve optimization aspects would be also very helpful by exploring the characterization of the PSO family with respect to typical needs raised in optimizations problems.
%in a variety of application fields

\section{Acknowledgements}
This work has been partially supported by the European Project SmartC2Net (n. ICT--318023) and the TENACE PRIN Project (n. 20103P34XC) funded by the Italian Ministry of Education, University and Research.

\appendix
\section{Appendix}
\label{appendix}

We use equation~\eqref{eq:muYB} for deriving formula~\eqref{eq:muYBv} for $\mu_{V}(k)$ as follow 
\begin{align*} 
  \mu_{V}(k)&=\mu_{Y|B}(k)-\mu_{Y|B}(k-1) \\
  &=\frac{ x_{pi}(k-1)+ \delta_{L|X}(k-1) l_{pi} + \delta_{B|Y}(k-1) b_i}{1 + \delta_{L|X}(k-1) + \delta_{B|Y}(k-1)} - x_{pi}(k-1) \\
  &=\frac{ x_{pi}(k-1)+ \delta_{L|X}(k-1) l_{pi} + \delta_{B|Y}(k-1) b_i}{1 + \delta_{L|X}(k-1) + \delta_{B|Y}(k-1)} \\
  &\phantom{=}-\frac{x_{pi}(k-1) + \delta_{L|X}(k-1) x_{pi}(k-1) + \delta_{B|Y}(k-1) x_{pi}(k-1)}{1 + \delta_{L|X}(k-1) + \delta_{B|Y}(k-1)} \\
  &=\frac{\delta_{L|X}(k-1) (l_{pi}(k-1) - x_{pi}(k-1)) + \delta_{B|Y}(k-1) (b_i(k-1) - x_{pi}(k-1))}{1 + \delta_{L|X}(k-1) + \delta_{B|Y}(k-1)}.
\end{align*}

Formula~\eqref{eq:sigma2YBv} for $\sigma_{V}^2(k)$ is derived from~\eqref{eq:sigma2YB} as follow 
\begin{align*} 
  \sigma_{V}^2(k)&=\sigma_{Y|B}^2(k)+\sigma_{Y|B}^2(k-1) \nonumber \\
  &=\frac{\sigma_{pi}^2(k-1)}{1 + \delta_{L|X}(k-1) + \delta_{B|Y}(k-1)} + \sigma_{pi}^2(k-1) \\
  &=\frac{(2+\delta_{L|X}(k-1) + \delta_{B|Y}(k-1))\sigma_{pi}^2(k-1)}{1+\delta_{L|X}(k-1) + \delta_{B|Y}(k-1)}.
\end{align*} 

Formula~\eqref{eq:muYBv3c} for $\mu_{V}(k)$ is derived from equations~\eqref{eq:muYBv} and~\eqref{eq:muYB} as follow.

From~\eqref{eq:sigma2YB} we get
\begin{align}
\frac{1}{\sigma_{Y|B}^2(k)}&= \frac{1}{\sigma_{Y|B}^2(k-1)}+\frac{1}{\sigma_{L|X}^2}+\frac{1}{\sigma_{B|Y}^2}, \label{eq:1sigma2YBk} \\
\frac{1}{\sigma_{Y|B}^2(k-1)}&= \frac{1}{\sigma_{Y|B}^2(k)}-\frac{1}{\sigma_{L|X}^2}-\frac{1}{\sigma_{B|Y}^2}, \label{eq:1sigma2YBkm1} \\
\frac{\sigma_{Y|B}^2(k)}{\sigma_{Y|B}^2(k-1)}&= 1-\delta_{L|X}(k) + \delta_{B|Y}(k), \label{eq:1sigma2YBksigma2YBkm1} \\
\frac{\sigma_{Y|B}^2(k-1)}{\sigma_{Y|B}^2(k)}&= 1+\delta_{L|X}(k-1) + \delta_{B|Y}(k-1). \label{eq:1sigma2YBkm1sigma2YBk}
\end{align}

From~\eqref{eq:deltaLXBY} for each value of $k$ and $h$ we get
\begin{align}
  \sigma_{Y|B}^2(k)\delta_{L|X}(h)&=\frac{\sigma_{Y|B}^2(k) \sigma_{Y|B}^2(h)}{\sigma_{L|X}^2} =\delta_{L|X}(k)\sigma_{Y|B}^2(h) \label{eq:sigma2YBdeltaLX} \\
  \sigma_{Y|B}^2(k)\delta_{B|Y}(h)&=\frac{\sigma_{Y|B}^2(k) \sigma_{Y|B}^2(h)}{\sigma_{B|Y}^2} =\delta_{B|Y}(k)\sigma_{Y|B}^2(h). \label{eq:sigma2YBdeltaBY} 
\end{align}

From~\eqref{eq:muYBv} we obtain a new fomula for $\mu_{V}(k)$ as 
\begin{align} 
  \mu_{V}(k-1)&=\frac{1 - \delta_{L|X}(k-1) - \delta_{B|Y}(k-1)}{1 - \delta_{L|X}(k-1) - \delta_{B|Y}(k-1)} \mu_{V}(k-1) \nonumber\\
  &=\frac{1}{1 - \delta_{L|X}(k-1) - \delta_{B|Y}(k-1)} \mu_{V}(k-1) \nonumber\\
  &\phantom{=}-\frac{\delta_{L|X}(k-1) + \delta_{B|Y}(k-1)}{1 - \delta_{L|X}(k-1) - \delta_{B|Y}(k-1)} \mu_{V}(k-1)\nonumber\\
  &=\frac{\sigma_{Y|B}^2(k-1)}{\sigma_{Y|B}^2(k-2)}\frac{\sigma_{Y|B}^2(k-2)}{\sigma_{Y|B}^2(k-1)} \frac{1}{1 - \delta_{L|X}(k-1) - \delta_{B|Y}(k-1)} \mu_{V}(k-1) \nonumber\\
  &\phantom{=}-\frac{\delta_{L|X}(k-1) + \delta_{B|Y}(k-1)}{1 - \delta_{L|X}(k-1) - \delta_{B|Y}(k-1)} \mu_{V}(k-1)\nonumber\\
  &=\frac{\sigma_{Y|B}^2(k-1)}{\sigma_{Y|B}^2(k-2)} \frac{1 + \delta_{L|X}(k-2) + \delta_{B|Y}(k-2)}{1 - \delta_{L|X}(k-1) - \delta_{B|Y}(k-1)} \mu_{V}(k-1) \nonumber\\
  &\phantom{=}-\frac{\delta_{L|X}(k-1) + \delta_{B|Y}(k-1)}{1 - \delta_{L|X}(k-1) - \delta_{B|Y}(k-1)} \mu_{V}(k-1), \text{ from \eqref{eq:1sigma2YBkm1sigma2YBk}}\nonumber\\
  &=\frac{\sigma_{Y|B}^2(k-1)}{\sigma_{Y|B}^2(k-2)} \frac{1 + \delta_{L|X}(k-2) + \delta_{B|Y}(k-2)}{1 - \delta_{L|X}(k-1) - \delta_{B|Y}(k-1)} \nonumber\\
  &\phantom{=} \frac{\delta_{L|X}(k-2) (l_{pi}(k-2) - x_{pi}(k-2)) + \delta_{B|Y}(k-2) (b_i(k-2) - x_{pi}(k-2))}{1 + \delta_{L|X}(k-2) + \delta_{B|Y}(k-2)} \nonumber\\
  &\phantom{=}-\frac{\delta_{L|X}(k-1) + \delta_{B|Y}(k-1)}{1 - \delta_{L|X}(k-1) - \delta_{B|Y}(k-1)} (x_{pi}(k-1)-x_{pi}(k-2)), \text{ from~\eqref{eq:muYBv}} \nonumber\\
  &=\frac{\sigma_{Y|B}^2(k-1)}{\sigma_{Y|B}^2(k-2)} \nonumber\\
  &\phantom{=} \frac{\delta_{L|X}(k-2) (l_{pi}(k-2) - x_{pi}(k-2)) + \delta_{B|Y}(k-2) (b_i(k-2) - x_{pi}(k-2))}{1 - \delta_{L|X}(k-1) - \delta_{B|Y}(k-1)} \nonumber\\
&\phantom{=}+\frac{\delta_{L|X}(k-1) + \delta_{B|Y}(k-1)}{1 - \delta_{L|X}(k-1) - \delta_{B|Y}(k-1)} (x_{pi}(k-2)-x_{pi}(k-1)) \nonumber\\
  &=\frac{\sigma_{Y|B}^2(k-2)}{\sigma_{Y|B}^2(k-2)} \nonumber\\
  &\phantom{=} \frac{\delta_{L|X}(k-1) (l_{pi}(k-2) - x_{pi}(k-2)) + \delta_{B|Y}(k-2) (b_i(k-2) - x_{pi}(k-2))}{1 - \delta_{L|X}(k-1) - \delta_{B|Y}(k-1)} \nonumber\\
  &\phantom{=}+\frac{\delta_{L|X}(k-1) + \delta_{B|Y}(k-1)}{1 - \delta_{L|X}(k-1) - \delta_{B|Y}(k-1)} (x_{pi}(k-2)-x_{pi}(k-1)), \text{ from~\eqref{eq:sigma2YBdeltaLX} and~\eqref{eq:sigma2YBdeltaBY}} \nonumber\\
  &=\frac{\delta_{L|X}(k-1) l_{pi}(k-2) - \delta_{L|X}(k-1) x_{pi}(k-2)}{1 - \delta_{L|X}(k-1) - \delta_{B|Y}(k-1)} \nonumber\\
  &\phantom{=}+\frac{\delta_{B|Y}(k-2) b_i(k-2) - \delta_{B|Y}(k-2) x_{pi}(k-2)}{1 - \delta_{L|X}(k-1) - \delta_{B|Y}(k-1)} \nonumber\\
  &\phantom{=}+\frac{\delta_{L|X}(k-1)x_{pi}(k-2) + \delta_{B|Y}(k-1)x_{pi}(k-2)}{1 - \delta_{L|X}(k-1) - \delta_{B|Y}(k-1)} \nonumber\\
  &\phantom{=}+\frac{-\delta_{L|X}(k-1)x_{pi}(k-1) - \delta_{B|Y}(k-1)x_{pi}(k-1)}{1 - \delta_{L|X}(k-1) - \delta_{B|Y}(k-1)} \nonumber\\
  &=\frac{\delta_{L|X}(k-1) (l_{pi}(k-2) - x_{pi}(k-1))+\delta_{B|Y}(k-1)(b_i(k-2)-x_{pi}(k-1))}{1 - \delta_{L|X}(k-1) - \delta_{B|Y}(k-1)} \label{eq:muYBvl2x1}
\end{align}

Finally from~\eqref{eq:muYBv} and~\eqref{eq:muYBvl2x1} we derive formula~\eqref{eq:muYBv3c} for $\mu_{V}(k)$ as 
\begin{align} 
  \mu_{V}(k)&=\frac{\delta_{L|X}(k-1) (l_{pi}(k-1) - x_{pi}(k-1)) + \delta_{B|Y}(k-1) (b_i(k-1) - x_{pi}(k-1))}{1 + \delta_{L|X}(k-1) + \delta_{B|Y}(k-1)} \nonumber\\
  &=\frac{\delta_{L|X}(k-1) l_{pi}(k-1) - \delta_{L|X}(k-1) x_{pi}(k-1)}{1 + \delta_{L|X}(k-1) + \delta_{B|Y}(k-1)}  \nonumber\\
  &\phantom{=}+\frac{\delta_{B|Y}(k-1) b_i(k-1) - \delta_{B|Y}(k-1) x_{pi}(k-1)}{1 + \delta_{L|X}(k-1) + \delta_{B|Y}(k-1)} \nonumber\\
  &\phantom{=}+\frac{\delta_{L|X}(k-1) l_{pi}(k-2) -\delta_{L|X}(k-1) l_{pi}(k-2)}{1 + \delta_{L|X}(k-1) + \delta_{B|Y}(k-1)} \nonumber\\
  &\phantom{=}+\frac{\delta_{B|Y}(k-1) b_i(k-2) - \delta_{B|Y}(k-1) b_i(k-2)}{1 + \delta_{L|X}(k-1) + \delta_{B|Y}(k-1)} \nonumber\\
  &=\frac{1 - \delta_{L|X}(k-1) - \delta_{B|Y}(k-1)}{1 - \delta_{L|X}(k-1) - \delta_{B|Y}(k-1)} \nonumber\\
  &\phantom{=}\left(\frac{\delta_{L|X}(k-1) l_{pi}(k-2) - \delta_{L|X}(k-1) x_{pi}(k-1)}{1 + \delta_{L|X}(k-1) + \delta_{B|Y}(k-1)} \right. \nonumber\\
  &\phantom{=}\left. +\frac{\delta_{B|Y}(k-1) b_i(k-2) - \delta_{B|Y}(k-1) x_{pi}(k-1)}{1 + \delta_{L|X}(k-1) + \delta_{B|Y}(k-1)} \right)\nonumber\\
  &\phantom{=}+\frac{\delta_{L|X}(k-1) l_{pi}(k-1) -\delta_{L|X}(k-1) l_{pi}(k-2)}{1 + \delta_{L|X}(k-1) + \delta_{B|Y}(k-1)} \nonumber\\
  &\phantom{=}+\frac{\delta_{B|Y}(k-1) b_i(k-1) - \delta_{B|Y}(k-1) b_i(k-2)}{1 + \delta_{L|X}(k-1) + \delta_{B|Y}(k-1)} \nonumber\\
  &=\frac{1 - \delta_{L|X}(k-1) - \delta_{B|Y}(k-1)}{1 + \delta_{L|X}(k-1) + \delta_{B|Y}(k-1)} \nonumber\\
  &\phantom{=}\left(\frac{\delta_{L|X}(k-1) l_{pi}(k-2) - \delta_{L|X}(k-1) x_{pi}(k-1)}{1 - \delta_{L|X}(k-1) - \delta_{B|Y}(k-1)} \right. \nonumber\\
  &\phantom{=}\left. +\frac{\delta_{B|Y}(k-1) b_i(k-2) - \delta_{B|Y}(k-1) x_{pi}(k-1)}{1 - \delta_{L|X}(k-1) - \delta_{B|Y}(k-1)} \right)\nonumber\\
  &\phantom{=}+\frac{\delta_{L|X}(k-1) l_{pi}(k-1) -\delta_{L|X}(k-1) l_{pi}(k-2)}{1 + \delta_{L|X}(k-1) + \delta_{B|Y}(k-1)} \nonumber\\
  &\phantom{=}+\frac{\delta_{B|Y}(k-1) b_i(k-1) - \delta_{B|Y}(k-1) b_i(k-2)}{1 + \delta_{L|X}(k-1) + \delta_{B|Y}(k-1)} \nonumber\\
  &=\frac{1 - \delta_{L|X}(k-1) - \delta_{B|Y}(k-1)}{1 + \delta_{L|X}(k-1) + \delta_{B|Y}(k-1)} \nonumber\\
  &\phantom{=}\frac{\delta_{L|X}(k-1) (l_{pi}(k-2) - x_{pi}(k-1))+\delta_{B|Y}(k-1) (b_i(k-2) - x_{pi}(k-1))}{1 - \delta_{L|X}(k-1) - \delta_{B|Y}(k-1)} \nonumber\\
  &\phantom{=}+\frac{\delta_{L|X}(k-1) (l_{pi}(k-1) - l_{pi}(k-2))+\delta_{B|Y}(k-1) (b_i(k-1) - b_i(k-2))}{1 + \delta_{L|X}(k-1) + \delta_{B|Y}(k-1)} \nonumber\\
  &=\frac{1 - \delta_{L|X}(k-1) - \delta_{B|Y}(k-1)}{1 + \delta_{L|X}(k-1) + \delta_{B|Y}(k-1)} \mu_{V}(k-1)\nonumber\\
  &\phantom{=}+\frac{\delta_{L|X}(k-1) (l_{pi}(k-1) - l_{pi}(k-2))}{1 + \delta_{L|X}(k-1) + \delta_{B|Y}(k-1)} \nonumber\\
  &\phantom{=}+\frac{\delta_{B|Y}(k-1) (b_i(k-1) - b_i(k-2))}{1 + \delta_{L|X}(k-1) + \delta_{B|Y}(k-1)}, \text{ from~\eqref{eq:muYBvl2x1}}
\end{align}

Another interesting formula for $\mu_{V}(k)$, it is derived from \eqref{eq:muYBvl2x1} observing that
\begin{equation}  l_{pi}(k-1)-l_{pi}(k-2) = (x_{pi}(k-1) - x_{pi}(k-2))=+(l_{pi}(k-1) - x_{pi}(k-1)) = -(l_{pi}(k-2) - x_{pi}(k-2)), \label{eq:lklkm1} 
\end{equation}
\begin{align}
  b_i(k-1)-b_i(k-2) = (x_{pi}(k-1) - x_{pi}(k-2)) \nonumber\\
  &\phantom{=}+(b_i(k-1) - x_{pi}(k-1)) \nonumber\\
  &\phantom{=}-(b_i(k-2) - x_{pi}(k-2)). \label{eq:bkbkm1}
\end{align}

Replacing~\eqref{eq:lklkm1} and~\eqref{eq:bkbkm1} in equation~\eqref{eq:muYBvl2x1} we get
\begin{align} 
  \mu_{V}(k-1)&=\frac{1 - \delta_{L|X}(k-1) - \delta_{B|Y}(k-1)}{1 + \delta_{L|X}(k-1) + \delta_{B|Y}(k-1)} \mu_{V}(k-1) \nonumber\\
  &\phantom{=}+\frac{\delta_{L|X}(k-1) (l_{pi}(k-1) - l_{pi}(k-2))}{1 + \delta_{L|X}(k-1) + \delta_{B|Y}(k-1)} \nonumber\\
  &\phantom{=}+\frac{\delta_{B|Y}(k-1) (b_i(k-1) - b_i(k-2))}{1 + \delta_{L|X}(k-1) + \delta_{B|Y}(k-1)} \nonumber\\
  &=\frac{1 - \delta_{L|X}(k-1) - \delta_{B|Y}(k-1)}{1 + \delta_{L|X}(k-1) + \delta_{B|Y}(k-1)} \mu_{V}(k-1) \nonumber\\
  &\phantom{=}+\frac{\delta_{L|X}(k-1)}{1 + \delta_{L|X}(k-1) + \delta_{B|Y}(k-1)} (x_{pi}(k-1)-x_{pi}(k-2))\nonumber\\
  &\phantom{=}+\frac{\delta_{L|X}(k-1) ((l_{pi}(k-1) - x_{pi}(k-1)) - (l_{pi}(k-2)-x_{pi}(k-2))}{1 + \delta_{L|X}(k-1) + \delta_{B|Y}(k-1)} \nonumber\\
  &\phantom{=}+\frac{\delta_{B|Y}(k-1) }{1 + \delta_{L|X}(k-1) + \delta_{B|Y}(k-1)} (x_{pi}(k-1)-x_{pi}(k-2)) \nonumber\\
  &\phantom{=}+\frac{\delta_{B|Y}(k-1) ((b_i(k-1)-x_{pi}(k-1)) - (b_i(k-2)-x_{pi}(k-2)))}{1 + \delta_{L|X}(k-1) + \delta_{B|Y}(k-1)} \nonumber\\
  &=\frac{1 - \delta_{L|X}(k-1) - \delta_{B|Y}(k-1)}{1 + \delta_{L|X}(k-1) + \delta_{B|Y}(k-1)} \mu_{V}(k-1) \nonumber\\
  &\phantom{=}+\frac{\delta_{L|X}(k-1)}{1 + \delta_{L|X}(k-1) + \delta_{B|Y}(k-1)} \mu_{V}(k-1) \nonumber\\
  &\phantom{=}+\frac{\delta_{L|X}(k-1) ((l_{pi}(k-1) - x_{pi}(k-1)) - (l_{pi}(k-2)-x_{pi}(k-2))}{1 + \delta_{L|X}(k-1) + \delta_{B|Y}(k-1)} \nonumber\\
  &\phantom{=}+\frac{\delta_{B|Y}(k-1) }{1 + \delta_{L|X}(k-1) + \delta_{B|Y}(k-1)} \mu_{V}(k-1) \nonumber\\
  &\phantom{=}+\frac{\delta_{B|Y}(k-1) ((b_i(k-1)-x_{pi}(k-1)) - (b_i(k-2)-x_{pi}(k-2)))}{1 + \delta_{L|X}(k-1) + \delta_{B|Y}(k-1)} \nonumber\\
  &=\frac{1 - \delta_{L|X}(k-1) - \delta_{B|Y}(k-1)}{1 + \delta_{L|X}(k-1) + \delta_{B|Y}(k-1)} \mu_{V}(k-1) \nonumber\\
  &\phantom{=}+\frac{\delta_{L|X}(k-1)+\delta_{B|Y}(k-1)}{1 + \delta_{L|X}(k-1) + \delta_{B|Y}(k-1)} \mu_{V}(k-1) \nonumber\\
  &\phantom{=}+\frac{\delta_{L|X}(k-1) ((l_{pi}(k-1) - x_{pi}(k-1)) - (l_{pi}(k-2)-x_{pi}(k-2))}{1 + \delta_{L|X}(k-1) + \delta_{B|Y}(k-1)} \nonumber\\
  &\phantom{=}+\frac{\delta_{B|Y}(k-1) ((b_i(k-1)-x_{pi}(k-1)) - (b_i(k-2)-x_{pi}(k-2)))}{1 + \delta_{L|X}(k-1) + \delta_{B|Y}(k-1)} \nonumber\\
  &=\frac{1}{1 + \delta_{L|X}(k-1) + \delta_{B|Y}(k-1)} \mu_{V}(k-1) \nonumber\\
  &\phantom{=}+\frac{\delta_{L|X}(k-1) ((l_{pi}(k-1) - x_{pi}(k-1)) - (l_{pi}(k-2)-x_{pi}(k-2)))}{1 + \delta_{L|X}(k-1) + \delta_{B|Y}(k-1)} \nonumber\\
  &\phantom{=}+\frac{\delta_{B|Y}(k-1) ((b_i(k-1)-x_{pi}(k-1)) - (b_i(k-2)-x_{pi}(k-2)))}{1 + \delta_{L|X}(k-1) + \delta_{B|Y}(k-1)}. \label{eq:muYBvl1x1l2x2}
\end{align}

Equation~\eqref{eq:muYBvl1x1l2x2} shows also that, when the mean $\mu_{V}(k)$ of the velocity of a particle is greater than $0$, the particle changes direction with respect to the current velocity (moment) only when the distance between the current position and the personal or global best postition is changed with respect to the previous step. Otherwise, the magnitude of new velocity is obtained as a fraction of that of the current velocity.

\end{document}